\documentclass[letterpaper, 10 pt, conference]{ieeeconf}  

\IEEEoverridecommandlockouts                              
\newtheorem{thm}{Theorem}[section]

\newtheorem{lem}[thm]{Lemma}
\newtheorem{prop}[thm]{Proposition}
\newtheorem{definition}{Definition}[section]{\bf}{\it}
\newtheorem{asump}{Assumption}[section]

\newtheorem{example}{Example}[section]
\usepackage{graphics} 
\usepackage{epsfig} 
\usepackage{epstopdf}
\usepackage{times} 
\usepackage{amsmath} 
\usepackage{amssymb}  
\usepackage{amsfonts}  
\usepackage{subfig}
\usepackage{epstopdf}
\usepackage{graphicx} 
\usepackage{amssymb}
\usepackage{enumerate}
\usepackage{color}

\newcommand*{\QEDB}{\hfill\ensuremath{\square}}%

\title{\LARGE \bf
Inducing Uniform Asymptotic Stability in Non-Autonomous Accelerated Optimization Dynamics via Hybrid Regularization}

\author{Jorge I. Poveda and Na Li
\thanks{J. I. Poveda is with the Department of Electrical, Computer and Energy Engineering, at the University of Colorado, Boulder. Email: \textsl{jorge.poveda@colorado.edu}.}
\thanks{Na Li is with the School of Engineering and Applied Sciences at Harvard University. Email: \textsl{nali@seas.harvard.edu}.}
}

\begin{document}
\maketitle
\thispagestyle{empty}
\pagestyle{empty}
\begin{abstract}
There have been many recent efforts to study accelerated optimization algorithms from the perspective of dynamical systems. In this paper, we focus on the \emph{robustness} properties of the time-varying continuous-time version of these dynamics. These properties are critical for the implementation of accelerated algorithms in feedback-based control and optimization architectures. We show that a family of dynamics related to the continuous-time limit of Nesterov's accelerated gradient method can be rendered unstable under arbitrarily small bounded disturbances. Indeed, while solutions of these dynamics may converge to the set of optimizers, in general, this set may not be uniformly asymptotically stable. To induce uniformity, and robustness as a byproduct, we propose a framework where we regularize the dynamics by using resetting mechanisms that are modeled by well-posed hybrid dynamical systems. For these hybrid dynamics, we establish uniform asymptotic stability and robustness properties, as well as convergence rates that are similar to those of the non-hybrid dynamics. We finish by characterizing a family of discretization mechanisms that retain the main stability and robustness properties of the hybrid algorithms. 
\end{abstract}
\section{INTRODUCTION}
\label{sec_introduction}
In this paper we focus on robust and fast gradient-based algorithms for the optimization problem
\begin{equation}\label{opti_problem}
\min~~f(x),~~~x\in\mathbb{R}^n,
\end{equation}

\vspace{-0.1cm}
\noindent
where $f$ is a smooth convex function. This type of problems has received significant attention due to the variety of applications that require fast algorithms with scalable rates of convergence. Two well-known accelerated gradient methods are Nesterov gradient \cite{Nesterov_Book} and Heavy-ball \cite{Polyak64} methods, which have inspired many following work, e.g., \cite{Lessar_IQC}, \cite{DissipativityNesterov},  \cite{InexactNesterov1} to just name a few. While these results were initially developed for discrete-time systems, recent works have focused on the development of continuous-time  algorithms modeled as ordinary differential equations (ODEs), e.g., \cite{ODE_Nesterov}, \cite{Wibisono1}, \cite{Wilson18},  \cite{Krichene16}.  In particular, it was shown in \cite{ODE_Nesterov}, that the time-varying ODE
\begin{equation}\label{ODEN1}
\ddot{x}+\frac{\ell(p)}{t}\dot{x}+p^2ct^{p-2}\nabla f(x)=0
\end{equation}

\vspace{-0.1cm}
\noindent
can be seen as the limiting continuous-time system obtained from Nesterov's gradient method for the case when $\ell(p)=p+1$, $c=1/p^2$, and $p=2$. These results have been further generalized in \cite{Wibisono1,Wilson18} and \cite{Shi18}. Moreover, it was also recently shown in \cite{Jadbabaie_RK} that for the case when $\ell(p)=2p+1$ and $h(x)=0.5|x|^2$, Runge-Kutta discretization methods applied to \eqref{ODEN1} can generate discrete-time algorithms that achieve acceleration. While these results have been instrumental in the analysis and design of various optimization algorithms with provable acceleration and convergence properties, the study of the robustness properties of these algorithms has been considered only recently \cite{ROBUST_FAST_ODE,Faziyab18Siam,Mihailo19,InexactNesterov1,InexactNesterov2}. Indeed, as it has been noted in the literature, e.g., \cite{Candes_Restarting,Wilson18}, dynamics of the form \eqref{ODEN1} may become unstable under small disturbances or even under forward Euler discretization. Apart from the safety concerns that arise from implementing algorithms that can be rendered unstable by small disturbances, lack of structural robustness properties in optimization algorithms is problematic for feedback-based control and optimization architectures that use real-time corrupted measurements of states and/or gradients, e.g., real-time reinforcement learning \cite{RL_survey}, adaptive control \cite{Sastry:89}, model-free online optimization \cite{Poveda:16}, etc. On the other hand, the study of stability and robustness properties in dynamical systems is nontrivial. As it has been shown in \cite{Vidyasagar_Book} and \cite{Goebel:12}, continuous-time, discrete-time, and hybrid dynamical systems can generate trajectories that converge to a particular point, but which fail to render stable the same point under arbitrarily small disturbances, even in cases when the disturbances converge to zero exponentially fast \cite{Teel_hespanha_GES}.

Motivated by this background, in this paper we study the robustness properties of the accelerated ODE \eqref{ODEN1} under with respect to persistent disturbances in the states and dynamics. We show that, in general, system \eqref{ODEN1} lacks a strong convergence property, called \textit{uniform attractivity}, which has been typically used to certify robustness in time-varying dynamical systems \cite{ODE2}, \cite{Goebel:12}. In turn, lack of \textit{uniform} attractivity has been historically linked to potential lack of robustness to small disturbances, see for instance \cite{LoriaADAptive,Sastry:89}. In order to induce uniformity in the convergence, we propose to regularize the dynamics \eqref{ODEN1} by using the framework of hybrid systems \cite{Goebel:12}. The resulting regularized dynamics combine continuous-time and discrete-time dynamics, and can be seen as robust periodic and persistently non-periodic restarting mechanisms designed to induce not only convergence and acceleration, but also robust stability.  We note that while the idea of using restarting to improve the convergence performance of accelerated gradient dynamics has been studied in \cite{Candes_Restarting} and \cite{RestartingAlgo12} for discrete-time systems, in \cite{ODE_Nesterov} and \cite[Ch. 9]{Krichene16} for ODEs, and in \cite{HybridNesterov1} for hybrid systems, to the best of our knowledge there is a lack of study in the literature concerning accelerated gradient ODEs that render uniformly asymptotically stable the set of attractors, and for which strictly positive margins of robustness under arbitrarily small time-varying bounded disturbances can be established. Indeed, the results of this paper open the door for the development of mode-free accelerated optimization algorithms based on multi-time scale approximations \cite{PovedaNaLi_ES}. In addition, we show that, as a consequence of having robust stability properties, a family of \emph{regular} discretization mechanisms, which include forward Euler and k-th order Runge-Kutta methods as special cases, can be used to preserve the stability and robustness properties of the dynamics.

\subsection*{Notation}
Given a compact set $\mathcal{A}\subset\mathbb{R}^n$, and a column vector $x\in\mathbb{R}^n$, we define $|x|_{\mathcal{A}}:=\min_{y\in\mathcal{A}}|x-y|$. We use $\mathbb{B}$ to denote a closed unit ball of appropriate dimension, $\rho\mathbb{B}$ to denote a closed ball of radius $\rho>0$, and $\mathcal{X}+\rho\mathbb{B}$ to denote the union of all sets obtained by taking a closed ball of radius $\rho$ around each point in the set $\mathcal{X}$. The closure of a set $\mathcal{X}$ is denoted as $\overline{\mathcal{X}}$ and its convex hull is given by $\text{con}(\mathcal{X})$. A function $f:\mathbb{R}^n\to\mathbb{R}$ is said to be radially unbounded if $f(x)\to\infty$ as $|x|\to\infty$. A function $f$ is of class $\mathcal{F}_{L}$ if its gradient is globally Lipschitz continuous. A function $f$ is of class $\mathcal{F}_{\mu,L}$ if $f\in\mathcal{F}_L$ and $f$ is $\mu$-strongly convex. A function $\rho$ is of class $\mathcal{K}_{\infty}$ if it is continuous, zero at zero, strictly increasing, and satisfies $\rho(s)\to\infty$ as $s\to\infty$.
%

\section{ON THE UNIFORM CONVERGENCE PROPERTIES OF THE ACCELERATED GRADIENT DYNAMICS}
\label{sec_motivation}
\subsection{Nominal Accelerated Gradient Dynamics}
System \eqref{ODEN1} can be rewritten in the following state space representation with $x_1:=x$, and $x_2:=\dot{x}$: 
\begin{align}\label{generalODE10}
\dot{x}_1=x_2,~~~\dot{x}_2=-\frac{\ell(p)}{t}x_2-cp^2t^{p-2}\nabla f(x_1),
\end{align}
where $t\geq t_0>0$ and $\ell(p)>1$ for all $p$.  Alternatively, if one selects $x_1:=x$ and $x_2:=x+\frac{t}{\ell(p)-1}\dot{x}_1$, system \eqref{ODEN1} can be written as
\begin{equation}\label{generalODE20}
\dot{x}_1=\frac{\ell(p)-1}{t}\left(x_2-x_1\right),~\dot{x}_2=\frac{-cp^{2}t^{p-1}}{\ell(p)-1}\nabla f(x_1).
\end{equation}
%
Irrespective of the state space selection, under convexity and suitable smoothness assumptions on $f$, for $p\geq2$ solutions of \eqref{ODEN1} minimize the sub-optimality measure $\tilde{f}(x_1(t)):=f(x_1(t))-f(x^*)$ at a rate $\mathcal{O}(1/t^p)$ \cite{ODE_Nesterov,Wibisono1}. Under further conditions on $f$ it can also be established that $x(t)$ converges to $x^*$ \cite{ROBUST_FAST_ODE}. However, while this type of convergence results are instrumental for the understanding of system \eqref{generalODE10}, they do not provide information related to the robustness properties of the system under small but persistent time-varying disturbances, which could be of adversarial nature. 

\subsection{Perturbed Accelerated Gradient ODE}
\label{sec_perturbed_ODEs}
Consider the ODE $\dot{x}=F(t,x)$ with state space representation \eqref{generalODE10} or \eqref{generalODE20}, and let $e_s(t),e_a(t):\mathbb{R}_{\geq0}\to\mathbb{R}^{2n+1}$ be measurable perturbation functions satisfying $|e_s(t)|\leq \varepsilon$ and $|e_a(t)|\leq \varepsilon$  for all $t\geq0$, for some $\varepsilon>0$. The perturbed version of equations \eqref{generalODE10} and \eqref{generalODE20} is then given by
\begin{equation}\label{ODEperturbed_1}
\dot{x}=F(t,x+e_s)+e_a.
\end{equation}
Typically, the signal $e_s(t)$ is related to unavoidable measurement noise that emerges in practical applications. On the other hand, the signal $e_a(t)$ captures model uncertainty, error approximations on the gradients $\nabla f(t)$, or offsets in the ``clock'' $\tau$ that coordinates the dynamics. Since we only impose an upper bound on $|e(t)|$, the signals $e_s(t)$ and $e_a(t)$ could also be of \emph{adversarial} nature.

For perturbed systems of the form \eqref{ODEperturbed_1}, we are interested in establishing the existence of a positive margin of robustness $\varepsilon$ such that all solutions of \eqref{ODEperturbed_1} behave in a similar way to the nominal dynamics \eqref{ODEN1}. In order to achieve this, a traditional approach is to establish \emph{uniform stability} and \emph{uniform convergence} properties for the nominal system.
\begin{definition}\label{defUGS}\cite[Def. 4.4]{khalil:book}
For the system $\dot{x}=f(t,x)$ with $x\in\mathbb{R}^n$, the origin $x^*=0$ is said to be \emph{uniformly stable} (US) if for each $\varepsilon>0$ there exists a $\delta_{\varepsilon}>0$ (independent of $t_0$), such that all solutions satisfying $|x(t_0)|\leq \delta$ also satisfy $|x(t)|\leq \varepsilon,~\forall~t\geq t_0$. It is said to be uniformly globally stable (UGS) if $\delta_{\varepsilon}$ can be chosen to satisfy $\lim_{\varepsilon\to\infty}\delta_{\varepsilon}=\infty$. \QEDB
 \end{definition}
\begin{definition}\label{defUGA}\cite[Def. 4.4]{khalil:book}
For the system $\dot{x}=f(t,x)$ with $x\in\mathbb{R}^m$, the origin $x^*=0$ is said to be \emph{uniformly globally attractive} (UGA) if for each pair $r>\varepsilon>0$ there exists a $T>0$ such that all solutions satisfying $|x(0)|\leq r$ also satisfy $|x(t)|\leq\varepsilon,~\forall~t\geq t_0+T$. \QEDB
\end{definition}

In words, the property of UGA asks that, for each pair $r>\varepsilon$, all solutions $x$ initialized in a $r$-neighborhood of the origin, must converge to $\varepsilon$-neighborhoods of the origin before some finite time $T+t_0$, with $T$ depending only on $(r,\varepsilon)$. It should be noted that UGA is a stronger notion compared to the classic notions of convergence used in optimization algorithms which usually do not impose any condition on how the convergence depends on the initial data $t_0$ and $z(0)$.
\begin{definition}\label{defUGAS}\cite[Def. 4.4]{khalil:book}
For the system $\dot{x}=f(t,x)$, the origin $x^*=0$ is said to be \emph{uniformly globally asymptotically stable} (UGAS) if it is uniformly globally stable and uniformly globally attractive. \QEDB
\end{definition}

The property of UGAS is relevant for the study of the robustness properties of dynamical systems, see for instance \cite{ODE2,Sastry:89,RelaxedPE}. Indeed, UGAS ensures robustness properties via the existence of converse Lyapunov functions. Unfortunately, as the following counter example shows, the accelerated gradient ODE \eqref{generalODE10} may generate trajectories that converge to the solution of \eqref{opti_problem} in a non-uniform way, even when $f$ is strongly convex.

\vspace{0.1cm}
\begin{example}\label{example1}
Consider the accelerated gradient ODE \eqref{generalODE10} in explicit time-varying form with $t_0\geq1$, $p=2$, $c=1$, and $f(x)=\frac{1}{2p^2}x_1^2$. Let $s=t-1$, and consider the time-varying dynamics in the $s$-time scale, given by
\begin{equation}\label{generalODE1c}
\begin{array}{l}
\dfrac{dx_1}{ds}=x_2\\
\dfrac{dx_2}{ds}=-h(s)x_2-p^2\nabla f =-h(s)x_2-x_1
\end{array},
\end{equation}
with $s\geq s_0\geq0$ and $h(s):=\frac{\ell(2)}{s+1}$. This system renders the origin $x^*=0$ UGS via the Lyapunov function $V(x)=0.5x_1^2+0.5x_2^2$, which leads to $\dot{V}(x)=-h(s)x_2^2\leq0$. Moreover, by the results of \cite{ODE_Nesterov} and \cite{Jadbabaie_RK}, all solutions satisfy $x(t)\to0$ as $t\to\infty$. To show that $\mathcal{A}$ is not UGAS, we can use the notion of limiting equations proposed by Artstein in \cite{ArtsteinLimiting}. Indeed, by \cite[Thm. 5.2]{ArtsteinLimiting}, if $h(s)\geq0$ and $\int h(s)ds$ is uniformly continuous in $\mathbb{R}_{\geq0}$, then system \eqref{generalODE1c} renders the origin UGAS if and only if there is no sequence $\{s_k\}^{\infty}_{k=1}$ with $s_k\to\infty$ as $k\to\infty$ such that $\int^{s_k+r}_{s_k}h(s)ds\to 0$ as $k\to\infty$ for each $r$. For system \eqref{generalODE1c}, $h(s)>0$ for all $s\geq0$, and since $h(s)\leq \ell(2)$ for all $s\geq0$, we have that $\int h(s)ds$ is uniformly continuous in $\mathbb{R}_{\geq0}$. However, for the function $h(s)$ in \eqref{generalODE1c} we have that $\int^{s_k+r}_{s_k}h(s)ds=\ell(2)\log\left(1+\frac{1}{s_k+r}\right)$, which converges to $0$ for all $r\geq0$ and any sequence satisfying $s_k\to\infty$. Therefore the origin is not UGAS for system \eqref{generalODE1c}.   \QEDB
\end{example}

\begin{figure}[t!]
\includegraphics[width=8.5cm,height=6.2cm]{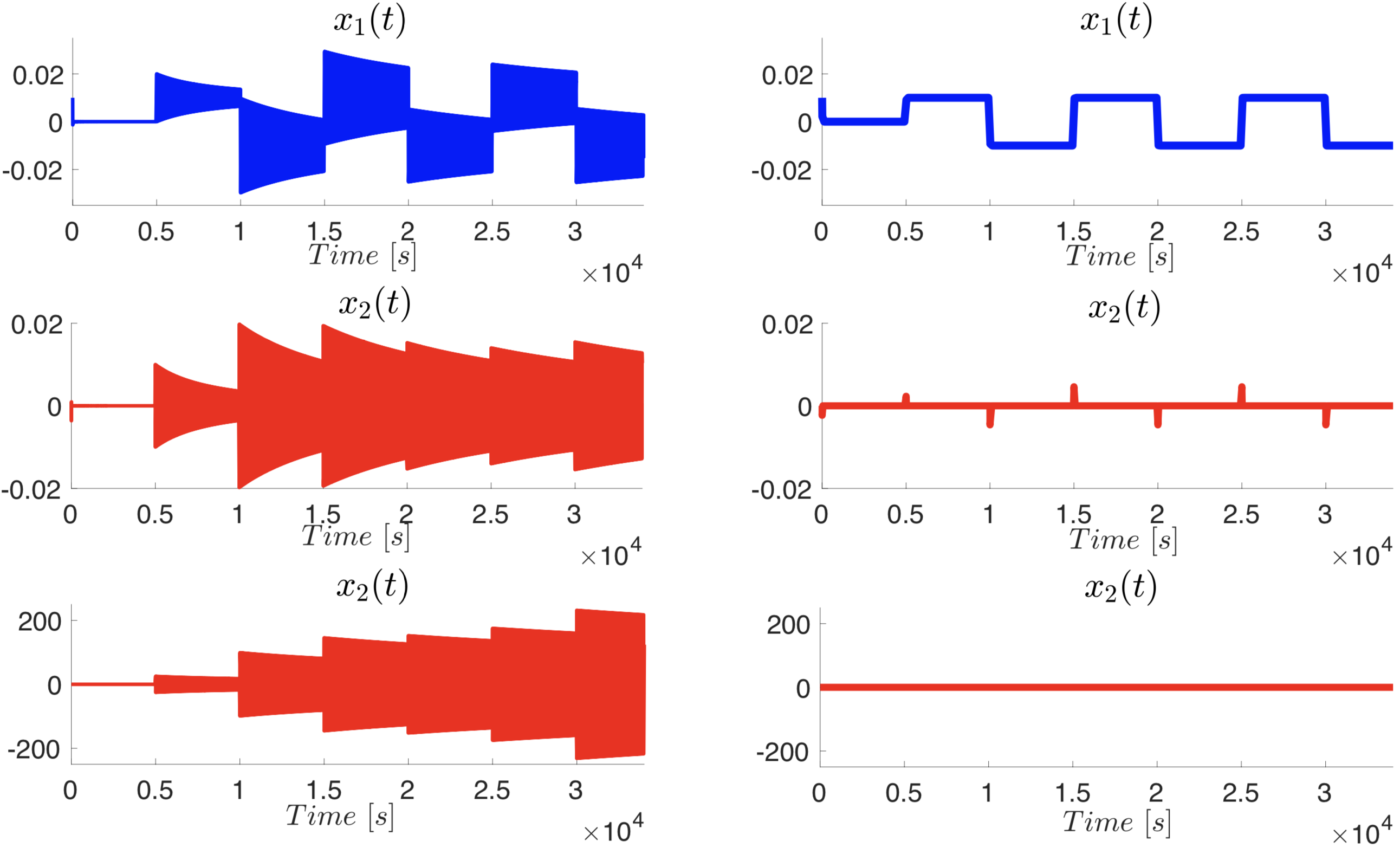}
\caption{Left plots: Instability in the dynamics \eqref{generalODE10} and \eqref{generalODE20} induced by disturbance $e_a(t)$. Right plots: Trajectories of \eqref{generalODE10} and \eqref{generalODE20} obtained when the hybrid regularization mechanism HAND-2 proposed in Section \ref{section2_convex} is implemented.}
\label{fig1}
\end{figure}

 Lack of uniformity with respect to the initial time in the convergence properties of the ODE \eqref{generalODE1c} implies that as $s\to\infty$, the damping term  $h(s)x_2$ takes longer and longer to react to small changes in the system. The left plots of Figure \ref{fig1} shows the trajectories of $x_1$ and $x_2$ generated under $\varepsilon$-disturbances on $\nabla f$ in \eqref{generalODE10} and \eqref{generalODE20}, with $\varepsilon=1\times10^{-3}$, $p=2$, $c=0.25$, $e_a(t)$ being a periodic square signal with period of $1\times10^4$ s, and $f(x)$ as in Example \ref{example1}. As shown in the plot, the $\varepsilon$-perturbation induces instability of the origin. Note that other type of arbitrarily small and state dependent adversarial disturbances $e_a$ and $e_s$ could also be considered. The robustness issues of equation \eqref{generalODE10} motivates us to study in the next section a class of regularization mechanisms that induce robust asymptotic stability properties in the accelerated ODEs. An example of the robust behavior induced by one these mechanisms is shown in the right plots of Figure \ref{fig1}, where we show the evolution of $x_1$ and $x_2$ under the same adversarial signal as in the left plots. 
%

\section{REGULARIZING THE GRADIENT DYNAMICS VIA ROBUST HYBRID SYSTEMS}
\label{sec_Nesterov}
In this section, we consider a class of regularization mechanisms that combine continuous-time and discrete-time dynamics. These mechanisms are modeled as hybrid dynamical systems (HDS) \cite{Goebel:12} with state $z\in\mathbb{R}^{2n+1}$ and dynamics

\vspace{-0.5cm}
\begin{subequations}\label{general_HDS}
\begin{align}
&z\in C,~~~~~\dot{z}=F(z),\label{flows_C}\\
&z\in D,~~~z^+=G(z)\label{jumps_D}.
\end{align}
\end{subequations}
System \eqref{general_HDS} exhibits continuous-time flows, given by \eqref{flows_C}, and discrete-time jumps, given by \eqref{jumps_D}. Therefore, their solutions $z$ are parameterized by a continuous-time index $t\in\mathbb{R}_{\geq0}$, which increases continuously during the flows, and a discrete-time index $j\in\mathbb{Z}_{\geq0}$, which increases by one during the jumps\footnote{For a precise definition of solutions for \eqref{general_HDS} we refer the reader to \cite{Goebel:12}.}. For HDS of the form \eqref{general_HDS}, an equivalent notion of UGAS can be defined for compact sets.
\begin{definition}\cite[Def. 3.6]{Goebel:12}
A compact set $\mathcal{A}\subset\mathbb{R}^{2n+1}$ is said to be UGAS if it satisfies the following two conditions:
\begin{enumerate}[(a)]
\item \textsl{Uniform global stability:} There exists a class $\mathcal{K}_{\infty}$ function $\alpha$ such that any solution $z$ to \eqref{general_HDS} satisfies $|z(t,j)|_{\mathcal{A}}\leq\alpha(|z(0,0)|_{\mathcal{A}})$ for all $(t,j)\in\text{dom}(z)$.
\item \textsl{Uniform global attractivity:} For each $\varepsilon>0$ and $r>0$ there exists $T>0$ such that, for any solution $z$ to \eqref{general_HDS} with $|z(0,0)|_{\mathcal{A}}\leq r$, $(t,j)\in\text{dom}(z)$ and $t+j\geq T$ imply $|z(t,j)|_{\mathcal{A}}\leq \varepsilon$. \QEDB
\end{enumerate}
\end{definition}

\vspace{0.1cm}
\noindent
Using the formalism \eqref{general_HDS}, and the state space representation \eqref{generalODE20}, we consider a family of regularized Hybrid Accelerated Nesterov Dynamics (HANDs) with overall state $z=[x^\top_1,x^\top_2,\tau]^\top\in\mathbb{R}^{2n+1}$ and hybrid dynamics 
\begin{subequations}\label{general_HANDS}
\begin{align}
&~~\dot{z}=F(z):=\left(\begin{array}{c}
\frac{2}{\tau}\left(x_2-x_1\right)\\
-2c\tau\nabla f(x_1),\\
1
\end{array}\right),~~z\in C,\label{flows_periodic1}\\
&z^+=G(z):=
\left(\begin{array}{l}
G_x(z)\\
G_{\tau}(z)
\end{array}\right),~~~~~~~~~~~~z\in D,
\end{align}
\end{subequations}
where for simplicity we used $\ell(p):=p+1$ and $p=2$, and where the mappings $G_x$ and $G_{\tau}$ are resetting functions to be designed. As in the non-hybrid case \eqref{generalODE10}, in order to study the robustness properties of \eqref{general_HANDS} we also consider perturbed HANDs of the form
\begin{subequations}\label{general_HANDS_perturbed}
\begin{align}
\dot{z}=F(z+e_1)+e_2,~~~~z+e_3\in C,\\
z^+=G(z+e_4)+e_5,~~~~z+e_6\in D,
\end{align}
\end{subequations}
where the signals $e_i:\mathbb{R}_{\geq0}\to\mathbb{R}^{2n+1}$ are all measurable admissible perturbations that satisfy $\sup_{t\geq0}|e_i(t)|\leq \varepsilon$ for some $\varepsilon>0$, for all $i\in\{1,2,\ldots,6\}$. By designing different types of flow and jump sets $C$ and $D$, as well as mappings $G_x$ and $G_{\tau}$, we will obtain hybrid systems generating periodic and aperiodic solutions for convex and strongly convex cost functions. 
%
%
%
%
\subsection{Hybrid Regularization for Radially Unbounded Convex Functions with Unique Minimizers}
We start by considering a regularization mechanism for cost functions $f$ satisfying the following assumption.
\begin{asump}\label{assumption_simple}
The function $f(x)$ is twice continuously differentiable, convex, radially unbounded, and has a unique minimizer $x^*\in\mathbb{R}^n$.\QEDB
\end{asump}

For cost functions satisfying Assumption \ref{assumption_simple}, we model the time index $\tau$ as a resetting clock, which leads to a HAND-1 \eqref{general_HANDS} with flows given by \eqref{flows_periodic1}, jumps $G(z)$ given by 
\begin{equation}\label{jumps_periodic1}
x^+= G_x(z)=x,~~~~\tau^+= G_{\tau}(z):=T_{\min},
\end{equation}
and flow and jump sets given by
\begin{subequations}\label{flow_jump_sets_1}
\begin{align}
C:=&\left\{z\in \mathbb{R}^{2n+1}:\tau\in [T_{\min},T_{\max}]\right\}\\
D:=&\left\{z \in \mathbb{R}^{2n+1}:\tau\in [T_{\text{med}},T_{\max}]\right\},
\end{align}
\end{subequations}
where $0<T_\text{min}<T_{\text{med}}\leq T_{\max}<\infty$. The resulting HDS \eqref{general_HANDS} is well-posed because the sets $C$ and $D$ are closed, the mappings $F$ and $G$ are continuous, $C\subset\text{dom}(F)$, and $D\subset\text{dom}(G)$. This follows by Thm. 6.30 in \cite{Goebel:12}.
%
%
%
Moreover, since $C\cap D\neq\emptyset$ the construction of the HAND-1 allows for non-unique solutions. In particular, the system allows resets of the clock at any instance such that the condition $\tau\geq T_{\text{med}}$ holds, but not later than when $\tau=T_{\max}$. For the particular case when the parameters are selected such that $T_{\text{med}}=T_{\text{max}}$, the resettings are periodic and the solutions generated by the HAND-1 are unique.  

It turns out that the simple modifications induced by \eqref{jumps_periodic1} and \eqref{flow_jump_sets_1} lead to a family of gradient algorithms that render UGAS the compact set
\begin{equation}\label{compact_set_1}
\mathcal{A}:=\{x^*\}\times\{x^*\}\times[T_{\min},T_{\max}],
\end{equation}
with strictly positive margins of robustness. The proof is presented in the Appendix.
\begin{thm}\label{thm1}
Suppose that Assumption \ref{assumption_simple} holds and consider the HAND-1. Then, the following holds:
\begin{enumerate}[(a)]
\item  Every maximal solution is complete and the set $\mathcal{A}$, given by \eqref{compact_set_1},  is UGAS.
\item For each $\delta\in\mathbb{R}_{>0}$ and each compact set $K_0\subset\mathbb{R}^{2n}$ there exists an $\varepsilon^*\in\mathbb{R}_{>0}$ and a $T\in\mathbb{R}_{>0}$ such that for every perturbation $e(t)$ satisfying $\sup_{t}|e(t)|\leq\varepsilon^*$ and every initial condition $z(0,0)\in K_0\times[T_{\min},T_{\max}]$ the solutions of the perturbed dynamics \eqref{general_HANDS_perturbed} satisfy $|z(t,j)|_{\mathcal{A}}\leq \delta$ for all $(t,j)\in\text{dom}(z)$ such that $t+j\geq T$.
\item For each $r\in\mathbb{R}_{>0}$ and each $z(0,0)$ such that $x_2(0,0)=x_1(0,0)$, $\tau(0,0)=T_{\min}$, and $x_1(0,0)\in K_0:=\{x^*\}+r\mathbb{B}$, we have that
\begin{align*}
f(x_1(t,0))-f^*\leq \frac{\beta}{t^2},
\end{align*}
for all $(t,j)\in\text{dom}(z)$ such that $j=0$, where $\beta:=\frac{r^2}{2c}+T^2_{\min}\tilde{f}(x_1(0,0))$ and $f^*=f(x^*)$.
 \QEDB
\end{enumerate}

\end{thm}

\vspace{0.2cm}
\noindent
In words, Theorem \ref{thm1} establishes that every solution generated by the HAND-1 will uniformly converge to the invariant compact set $\mathcal{A}$. Moreover, when $\varepsilon$-bounded perturbations of arbitrary frequency and/or adversarial nature are added to the states or dynamics of the system, the new solutions of the perturbed system will converge uniformly to the set $\mathcal{A}+\delta\mathbb{B}$, where $\delta>0$. To the knowledge of the authors, this type of robustness result has not been established before for the continuous-time accelerated gradient dynamics. Indeed, the proof of Theorem \ref{thm1} relies on an invariance principle for well-posed HDS that, to our knowledge, has not been used before in the analysis of accelerated optimization algorithms. Finally, item c) says that the sub-optimality measure $f(x_1)-f^*$ decreases at a rate of $\mathcal{O}(1/t^2)$ during the first interval of flow, which implies that, given $\delta>0$, if $T_{\text{med}}$ is selected such that
$$T_{\text{med}}\geq \sqrt{\frac{\beta}{\delta}}+T_{\min}>0,$$
then $f(x(t,j))-f^*\leq \delta,~~~\forall~(t,j)\in S_{T_{\text{med}}}$, where  $S_{T_{\text{med}}}:=\{(t,0)\in\text{dom}(z):t\geq T_{\text{med}}-T_{\text{min}}\}$.  However, as $T_{\text{med}}\to\infty$, the HAND-1 behaves as the time-varying ODE \eqref{generalODE10}, which establishes a clear tradeoff between acceleration and robustness. Note that while item (c) establishes acceleration only during the first interval of flow, it is possible to establish a semi-acceleration property for all $(t,j)\in\text{dom}(z)$ by generating similar bounds that hold for \emph{each} interval of flow (with different constants $\beta$), see \cite{PovedaNaLi_ES}.
\vspace{0.1cm}
\subsection{Hybrid Regularization for Strongly Convex Functions with Lipschitz Gradient}
\label{section2_convex}
We now consider cost functions $f(\cdot)$ that are also strongly convex and have a globally Lipschitz gradient.
\vspace{0.1cm}
\begin{asump}\label{assumption_HAND2}
The cost function $f$ is of class $\mathcal{F}_{\mu,L}$. \QEDB
\end{asump}
\vspace{0.2cm}

For functions of class $\mathcal{F}_{\mu,L}$ we are interested in designing HANDs with the UGAS property, and which, additionally, guarantee an exponential decay of the sub-optimality measure $f(x)-f^*$. To achieve this, we consider the continuous-time dynamics \eqref{flows_periodic1} combined with the following discrete-time dynamics:
\begin{subequations}\label{jumps_periodic2}
\begin{align}
x^+&= G_x(x,\tau):=[x^\top_1,x_1^\top]^\top,\\
\tau^+&= G(x,\tau):=T_{\min},
\end{align}
\end{subequations}
and flow and jump sets given by
\begin{subequations}\label{flow_jump_sets_2}
\begin{align}
C:=&\left\{z\in \mathbb{R}^{2n+1}:\tau\in [T_{\min},T_{\max}]\right\}\\
D:=&\left\{z \in \mathbb{R}^{2n+1}:\tau\in \{T_{\max}\}\right\},
\end{align}
\end{subequations}
where $0<T_\text{min}<T_{\max}<\infty$. Closedness of the sets $C$ and $D$, as well as continuity of the mappings $F$ and $G$ guarantee that the resulting HDS is also well-posed. Indeed, the HAND-2 describes an algorithm where the clock $\tau$ and the state $x_2$ are periodically reset to $T_{\min}$ and $x_1$, respectively, which is a typical resetting mechanism used in optimization algorithms with momentum \cite{Candes_Restarting}. The following theorem shows that this system also guarantees UGAS and robustness of the set \eqref{compact_set_1}, with an exponential decay in the sub-optimality measure, provided the jumps satisfy a quadratic dwell-time like condition. The proof is presented in the Appendix.
\begin{thm}\label{thm_HAND2}
Suppose that Assumptions \ref{assumption_simple} and \ref{assumption_HAND2} hold. Consider the HAND-2, and let $0<T_{\min}<T_{\max}<\infty$ and $c>0$ such that the following inequality is satisfied:
\begin{equation}\label{inequality_period2}
T^2_{\max}-T_{\min}^2>\frac{1}{\mu c}.
\end{equation}
Then, the following holds:
\begin{enumerate}[(a)]
\item Every maximal solution is complete, and the set $\mathcal{A}$, given by \eqref{compact_set_1}, is UGAS.
\item For each $\delta\in\mathbb{R}_{>0}$ and each compact set $K_0\subset\mathbb{R}^{2n}$ there exists $\varepsilon^*, T\in\mathbb{R}_{>0}$ such that for every perturbation $e(t)$ satisfying $\sup_{t}|e(t)|\leq\varepsilon^*$ and every initial condition $x(0,0)\in K_0$ the solutions of the perturbed hybrid dynamics \eqref{general_HANDS_perturbed} satisfy $|z(t,j)|_{\mathcal{A}}\leq \delta$ for all $(t,j)\in\text{dom}(z)$ such that $t+j\geq T$.
\item Let $\tilde{x}_{10}=x_1(0,0)-x^*$ and $\Delta T:=T_{\max}-T_{\min}$. If $x_1(0,0)=x_2(0,0)$ and $\tau(0,0)=T_{\min}$, the sub-optimality measure satisfies 
\begin{align}\label{sub_optimal_measure}
f(x_1(t,j))-f^*&\leq k_a \exp\left(-\tilde{k}_b\tilde{\alpha}(t+j) \right)  |\tilde{x}_1(0,0)|^2,
\end{align}
for all $(t,j)\in\text{dom}(z)$, where $k_a>0$, $\tilde{k}_b:=1-k_0$, 
\begin{equation}\label{constant_decay1}
k_0:=\frac{(c\mu)^{-1}+T_{\min}^2}{T^2_{\max}},
\end{equation}
and $\tilde{\alpha}(t+j):=\frac{\max\left\{t+j-\Delta T,0\right\}}{{\Delta T+1}}$.\QEDB
\end{enumerate}
\end{thm}

\vspace{0.2cm}
\noindent
Theorem \ref{thm_HAND2} states that condition \eqref{inequality_period2} is sufficient to guarantee UGAS and exponential decay of the sub-optimality measure. Indeed, this condition can be equivalently written as $T_{\min}+T_{\max}>\frac{1}{c\mu \Delta T}$. When $T_{\min}+T_{\max}>1$, this condition is satisfied if the following dwell-time condition holds
$$T_{\max}-T_{\min}>(c\mu)^{-1}.$$
%
%
%
%
%

\noindent
For cost functions satisfying Assumption 2 it is possible to establish additional properties for the HAND-2, e.g., uniform global \emph{exponential} stability. Such results are omitted in this paper and can be found in \cite{PovedaNaLi_ES}. 

For resetting mechanisms such as \eqref{jumps_periodic2}, it is useful to characterize the optimal switching frequency $\Delta T$ that minimizes the bound in \eqref{sub_optimal_measure} for a given window of time, see \cite{Candes_Restarting} for a discrete-time version of this result. In order to make this question tractable in our setting, we replace the constant $k_0$ in \eqref{constant_decay1} by the constant
\begin{equation}\label{k_1}
k_1:=\frac{(c\mu)^{-1}+T_{\min}^2}{\Delta T^2},
\end{equation}
which satisfies $k_1>k_0$, and we replace \eqref{inequality_period2} by the stronger condition
\begin{equation}\label{stronger_condition}
\Delta T^2-T_{\min}^2>\frac{1}{c\mu}.
\end{equation}
It is easy to see that condition \eqref{stronger_condition} implies condition \eqref{inequality_period2}. The proof of the following Lemma is also presented in the Appendix.
\begin{lem}\label{lemma_optimal_switching}
Let $k_1$ be given by \eqref{k_1} and suppose that condition \eqref{stronger_condition} holds. Let $(t,j)\in\text{dom}(z)$ be such that $t=j\Delta T$. Then, the optimal switching frequency $\Delta T^*$ that minimizes $k^{j}_1\big|_{j=t/\Delta T}$ is given by $\Delta T^*=e~\sqrt{\frac{1}{c\mu}+T^2_{\min}}$, and for each $\varepsilon>0$, this switching frequency guarantees that $f(x_1(t,j))-f^* \leq \varepsilon,~~\forall~(t,j)\in\text{dom}(z)$ such that $t=j\Delta T$ and $t\geq \frac{e}{2}\sqrt{\frac{1}{c\mu}+T^2_{\min}}\log\left(\frac{\tilde{f}(x_1(0,0))}{\varepsilon}\right)$. \QEDB
%
%
\end{lem}
\vspace{0.1cm}

Lemma \ref{lemma_optimal_switching} says that for any precision $\varepsilon>0$, the convergence time of the sub-optimality measure is of the order $O\left(\sqrt{\frac{1}{c\mu}+T^2_{\min}} \log\left({\frac{1}{\varepsilon}}\right)\right)$. A similar result is derived in \cite{Candes_Restarting} for the classic discrete-time Nesterov dynamics, and in \cite[Ch. 9]{Krichene16} and \cite{CDC19Ochoa} for continuous-time mirror descent.
%
%

\section{A STABLE DISCRETIZATION OF THE HYBRID DYNAMICS}
\label{sec_discretization}
%
In this section we show that, unlike their non-hybrid counterparts, the HANDs developed in Section \ref{sec_Nesterov} retain their (semi) global stability \emph{and} robustness properties under a variety of discretization mechanisms, including forward-Euler and k-Order Runge-Kutta integration schemes. This is in contrast to the unstable behavior that may emerge under simple Euler discretization of the non-hybrid ODE \eqref{ODEN1}. In our case, the stability properties of the discretized algorithms are inherited from the stability properties of the HANDs.

To model the discretized hybrid dynamics we use the framework of Hybrid Simulators \cite{SimulatorHybridSystems}, where a discretized hybrid system $\mathcal{H}_h$ is represented by the dynamics
\begin{equation}\label{discretized_HANDS}
z^+=F_h(z),~z\in C_h,~~\text{and}~~z^+=G_h(z),~z\in D_h,
\end{equation}
where the elements $(F_h,C_h,G_h,D_h)$ are obtained via a discretization mechanism with step size $h>0$. Unlike solutions of the HANDs considered in the previous section, which were defined on hybrid time domains, the solutions of \eqref{discretized_HANDS} are defined on discrete time domains\footnote{We refer the reader to \cite{SimulatorHybridSystems} for a complete description of hybrid simulators and their definition of solutions.}. To obtain ``well-posed'' discretized dynamics \eqref{discretized_HANDS}, we will consider a class of \emph{regular} discretization mechanisms.
\begin{definition}\label{regular_discretization}
The discretized HAND $\mathcal{H}_h$ is said to be \emph{regular} if the data $(F_h,C_h,G_h,D_h)$ satisfies the following conditions:
\begin{itemize}
\item  $F_h$ is such that, for each compact set $K\subset\mathbb{R}^n$, there exists a function $\rho\in\mathcal{K}_\infty$ and $h^*>0$ such that for each $z\in C_h\cap K$ and each $h\in(0,h^*]$
\begin{equation}\label{condition_discrete_flows}
F_h(z)\subset z+h~\overline{\text{con}} F(z+\rho(h)\mathbb{B})+h\rho(h)\mathbb{B}.
\end{equation}
\item $G_h$ is such that for any decreasing sequence $h_i\to0$ we have that $G_0=G(z)$, where $G_0$ is the graphical limit of $G_{h_i}$.
\item The sets $C_h$ and $D_h$ are such that for any positive monotone decreasing sequence $\{h_i\}^{\infty}_{i=1}$ such that $h_i\to 0$ we have that $\lim~\sup_{i\to\infty}~C_{h_i}\subset C$ and $\lim~\sup_{i\to\infty}~D_{h_i}\subset D$.  
\end{itemize}
\end{definition}
Examples of mappings $F_h$ satisfying the conditions of Definition \ref{regular_discretization} include forward-Euler and the consistent $S-$Order Runge-Kutta methods \cite{SimulatorHybridSystems}, given by $F_h(z)=z+hF(z)$, and
\begin{equation}\label{discretization2}
F_h(z)=z+h\sum_{k=1}^S b_kF(g_k),~~g_k=z+h\sum_{\ell=1}^{i-1}a_{ij}F(g_j),
\end{equation}
respectively, where $\sum_{k=1}^S b_k=1$, and $S=\{1,2,\ldots,\bar{s}\}$, $\bar{s}\in\mathbb{Z}_{>1}$. In order to obtain a discretized system whose updates do not abandon prematurely the set $C$ after a discretized flow, we can consider the discretized jump set $D_h=D\cup\{z:y\in C,z=F_h(y)\notin C\}$. Based on this we can consider the following Runge-Kutta-based hybrid discretizations for the HANDs studied in Section \ref{sec_Nesterov} that use $D_h$, $F_h$ given by \eqref{discretization2}, $C_h=C$, and $G_h=G$. Since regular discretization mechanisms with sufficiently small discretization steps preserve the convergence properties of hybrid systems (in a semi-global practical sense) we can obtain the next result, which follows directly by using items (a) and (b) of Theorems \ref{thm1} and \ref{thm_HAND2}, and applying \cite[Thm. 5.3]{SimulatorHybridSystems} to the HANDs 1 and 2.
\begin{prop}\label{theorem_discretization}
Consider the HANDs $\mathcal{H}$ of Section \ref{sec_Nesterov} under their corresponding assumptions. Then, for each $r>\varepsilon>0$ there exists a $h^*>0$ such that for all $h\in(0,h^*)$ there exists a $T>0$ such that if $|z_h(0,0)|_{\mathcal{A}}\leq r$ then $|z_h(k,j)|_{\mathcal{A}}\leq \varepsilon$ for all $(k,j)\in\text{dom}(z_h)$ such that $kh+j\geq T$. \QEDB
\end{prop}
\vspace{0.2cm}

The semi-global practical result of Proposition \ref{theorem_discretization} gives the existence of a sufficiently small upper bound $h^*>0$ for the step size, such that for any $h\in(0,h^*]$, the stability and uniform convergence properties of the HANDs are retained from compact sets of initial conditions. However, Theorem \ref{theorem_discretization} does not provide any constructive information for the selection of a step size $h$ that induces acceleration in discrete time. Related results in this direction have been recently pursued in \cite{Jadbabaie_RK}. 

%
 %
%
\section{Conclusions}
\label{sec_conclusions}
We studied robustness and uniform asymptotic stability properties of a class of \emph{time-varying} gradient ODEs related to the continuous-time limit of the Nesterov's algorithm. We showed that, even for strongly convex functions, the time-varying ODE may not render the set of optimizers UGAS, a property that traditionally has been used to certify robustness properties in feedback control systems. In order to induce this property in optimization dynamics, we propose two different regularization mechanisms based on well-posed hybrid dynamical systems, and we characterized the stability, convergence, and robustness properties of the emerging algorithms. To the knowledge of the authors these are the first results that establish robust UGAS for the time-varying Nesterov's ODE with momentum, preserving semi-acceleration properties for non-strongly convex functions. Finally, we showed that the UGAS property implies that a family of regular discretization mechanisms preserves the main properties of the hybrid dynamics for sufficiently small step size. 

\bibliographystyle{ieeetr}
\bibliography{Biblio}

\section{Proofs}
\label{sec_proofs}
\subsection{Proof of Theorem \ref{thm1}}
\label{proof_thm1}
Since $\tau$ is always constrained to evolve in the set $[T_{\min},T_{\max}]$, we have that $|z|_{\mathcal{A}}=|x-\textbf{1}_2\otimes x^*|$. Based on this, consider the following Lyapunov-like function studied also in \cite{Wibisono1,Wilson18}:
\begin{equation}\label{Lyapunov_function1}
V(z)=\frac{|x_2-x^*|^2}{2}+c\tau^2 (f(x_1)-f(x^*)),
\end{equation}
which, under Assumption \ref{assumption_simple}, is positive definite with respect to $\mathcal{A}$ and radially unbounded. Thus, there exists $\alpha_1,\alpha_2\in \mathcal{K}_{\infty}$ such that  $\alpha_1(|z|_{\mathcal{A}})\leq V(x)\leq \alpha_2(|z|_{\mathcal{A}})$ for all $z\in C\cup D$. The derivative of $V$ with respect to time satisfies
\begin{subequations}\label{inequality_flows_1}
\begin{align}
\dot{V}(z)=&\nabla V(z)^\top \dot{z}\notag\\
&=\left[c\tau^2\nabla f(x_1),~(x_2-x^*),~2c\tau (f(x_1)-f(x^*))\right]^\top\dot{z},\notag\\
&=2c\tau\nabla f(x_1)^\top(x_2-x_1)-2c\tau (x_2-x^*)^\top\nabla f(x_1)\notag\\
&~~~~+2c\tau (f(x_1)-f(x^*)),\notag\\
&=-2c\tau\left[\nabla f(x_1)^\top (x_1-x_1^*)- (f(x_1)-f(x^*))\right],\notag\\
&=:u_C(z)\leq 0,~~~~\forall~z\in C, 
\end{align}
\end{subequations}
where the last inequality follows by Assumption \ref{assumption_simple} which implies that $f(x^*)-f(x_1)-\nabla f(x_1)^\top (x^*-x_1)\geq0$. To show that this inequality is strict for all $x_1\neq x^*$, suppose by contradiction that there exists $x_1\neq x^*$ such that $f(x^*)-f(x_1)-\nabla f(x_1)^\top (x^*-x_1)=0$. Let $\alpha_1:=f(x_1)$ and define the set $\Omega_{\alpha_1}:=\{x\in\mathbb{R}^n: f(x)\leq \alpha_1\}$. Since any $x^*$ is optimal, we have that $f(x^*)=f^*\leq f(x_1)$ and therefore $x^*\in  \Omega_{\alpha_1}$. Since $f$ is twice continuously differentiable, the Hessian of $f$ is continuous and uniformly bounded on compact sets. Therefore, there exists $L_{\alpha_1}>0$ such that $|\nabla f(x'_1)-\nabla f(x''_1)|\leq L_{\alpha_1} |x'_1-x''_1|$ for all $(x'_1,x''_1)\in \Omega_{\alpha_1}$. By the convexity and the Lipschitz properties in $\Omega_{\alpha_1}$, we obtain
\begin{equation}
f(x^*)-f(x_1)-\nabla f(x_1)^{\top}(x^*-x_1)\geq \frac{1}{2L_{\alpha_1}}|\nabla f(x_1)|^2,
\end{equation}
but since by assumption the left hand side of the inequality is zero, we must have that $|\nabla f(x_1)|=0$, which is a contradiction given that $x_1\neq x^*$. Therefore $u_C(0)^{-1}=\{x^*\}\times\mathbb{R}^n\times[T_{\min},T_{\max}]$. 

On the other hand, the change of the Lyapunov-like function \eqref{Lyapunov_function1} during the jumps is given by
\begin{subequations}\label{inequality_jumps_1}
\begin{align}
V(z^+)-V(z)&=\frac{|x_2-x^*|^2}{2}+cT_{\min}^2 (f(x_1)-f(x^*))\notag\\
&~~~-\frac{|x_2-x^*|^2}{2}-c\tau^2 (f(x_1)-f(x^*)),\notag\\
&=-c(f(x_1)-f(x^*))(\tau^2 -T_{\min}^2),\notag\\
&=u_D(z)\leq 0,~~~~~\forall~z\in D,
\end{align}
\end{subequations}
where the last inequality follows by the fact that $\tau^2\geq T_{\text{med}}^2> T^2_{\min}$ in the set $D$. Therefore, the Lyapunov-like function \eqref{Lyapunov_function1} does not increase during jumps. Since the system is well-posed, the hybrid invariance principle \cite[Ch. 8]{Goebel:12} can be applied. Indeed, note that
$u^{-1}_C(0)=\{x^*\}\times\mathbb{R}^n\times[T_{\text{min}},T_{\text{max}}]$, $u^{-1}_D(0)=\left(\{x^*\}\times\mathbb{R}^n\times[T_{\text{min}},T_{\text{max}}]\right)\cup\left(\mathbb{R}^n\times\mathbb{R}^n\times\{T_{\text{min}}\}\right)$, and $G(u^{-1}_D(0))=\{x^*\}\times\mathbb{R}^n\times\{T_{\text{min}}\}$, where $G(z):=G_x(z)\times G_{\tau}(z)$. Let $U=C\cup D$, $r\in V(U)$ and 
\begin{equation}
\mathcal{W}:=V^{-1}(r)\cap U\cap [u^{-1}_C(0)\cup \left(u^{-1}_D(0)\cap G(u_D^{-1}(0))\right)].
\end{equation}
Then,
\begin{align}\label{subsetset}
\mathcal{W}_r=&\Big\{[x_1,x_2,\tau]\in\mathbb{R}^{2n+1}: x_1=x^*,~|x_2-x^*|=\sqrt{2r},\notag\\
&~~\tau\in[T_{\text{min}},T_{\text{max}}]\Big\},
\end{align}
and by \cite[Corollary 8.4]{Goebel:12} every solution of the HAND-1 approaches the largest weakly invariant subset of \eqref{subsetset}. To show that this set corresponds to the case $r=0$, note that by Assumption \ref{assumption_simple} any invariant solution in \eqref{subsetset} satisfies $x_1=x^*$ and $\dot{x}_1=0=2\tau^{-1}(x_2-x^*)$, which can only happen if $x_2=x^*$. Since the jumps satisfy $x^+=x$, the largest weakly invariant subset of \eqref{subsetset} corresponds $\mathcal{W}_0=\mathcal{A}$. Therefore, since $G(D)\subset C\cup D$, by the Barbasin-Krasovskii-LaSalle theorem for hybrid systems \cite[Thm. 8.8]{Goebel:12}, the set $\mathcal{A}$ is UGAS. Item (b) follows now directly by \cite[Lemma 7.20]{Goebel:12}. 

To show item (c), let $\tilde{f}(x_1)=f(x_1)-f^*$ and note that since the Lyapunov-like function \eqref{Lyapunov_function1} does not increase during flows or jumps, we have that $V(z(t_0+t,j))\leq V(z(t_0,j))$, for all $t\geq0$ such that $(t+t_0,j)\in\text{dom}(z)$ and $(t_0,j)\in\text{dom}(z)$, i.e., during flows. Taking $t=0$ and $j=0$ we get
\begin{align}
\tilde{f}(x_1(t,0))&\leq \frac{|\tilde{x}_1(0,0)|^2}{2c\tau^2}+T^2_{\min}\frac{\tilde{f}(x_1(0,0))}{\tau^2}\label{bound_thm1f}\\
&\leq \frac{r^2}{2c\tau^2}+T^2_{\min}\frac{\tilde{f}(x_1(0,0))}{\tau^2}\leq \frac{\beta}{\tau^2},
\end{align}
with $\beta:=\frac{r^2}{2c}+T^2_{\min}\tilde{f}(x_1(0,0))$. Thus, for any $\delta>0$ the condition $\tau>\sqrt{\beta/\delta}$ implies that $\tilde{f}(x_1(t,0))\leq \delta$, which can always be induced by designing $T_{\text{med}}$ and $T_{\max}$ such that $T_{\text{med}}-T_{\text{min}}\geq\sqrt{\beta/\delta}$. \null \hfill \null $\blacksquare$
%

\subsection{Proofs of Theorem \ref{thm_HAND2}} 
\label{proof_thm2}
To prove item (a), consider again the Lyapunov-like function \eqref{Lyapunov_function1}. Since the flow map is still given by \eqref{flows_periodic1}, and the cost function $f(\cdot)$ is strongly convex, inequality \eqref{inequality_flows_1} still holds. On the other hand, during jumps we now have
\begin{align}
V(z^+)-V(z)&=\frac{|x_1-x^*|^2}{2}+T_{\min}^2c \left(f(x_1)-f^*\right)\\
&~~~~-\frac{|x_2-x^*|^2}{2}-\tau^2 c\left(f(x_1)-f^*\right),\notag\\
&\leq-c\left(f(x_1)-f(x^*)\right)\left[\tau^2-T^2_{\min}-\frac{1}{\mu c}\right]\notag\\
&~~~-\frac{|x_2-x^*|^2}{2},\notag\\
&=:u_D(z)\leq 0,~~\forall~z\in D,
\label{inequality_jumps_4}
\end{align}
%
where we used the strong convexity of $f$ and inequality \eqref{inequality_period2} to get inequality \eqref{inequality_jumps_4}. Using again the fact that the system is nominally well-posed and the hybrid invariance principle of \cite[Ch. 8]{Goebel:12}, we now have $u^{-1}_C(0)=\{x^*\}\times\mathbb{R}^n\times[T_{\text{min}},T_{\text{max}}]$, $u^{-1}_D(0)=\left(\{x^*\}\times\{x^*\}\times[T_{\text{min}},T_{\text{max}}]\right)\cup\left(\mathbb{R}^n\times\{x^*\}\times\{(T^2_{\text{min}}+(c\mu)^{-1}))^{0.5}\}\right)$, and $G(u^{-1}_D(0))=\{x^*\}\times\{x^*\}\times\{T_{\text{min}}\}$, where $G(z):=G_x(z)\times G_{\tau}(z)$. Let $U=C\cup D$, $r\in V(U)$ and note that the set
\begin{equation}
\mathcal{W}:=V^{-1}(r)\cap U\cap [u^{-1}_C(0)\cup \left(u^{-1}_D(0)\cap G(u_D^{-1}(0))\right)]
\end{equation}
 is again given by \eqref{subsetset}. Since the largest weakly invariant set in $\mathcal{W}$ is again given by $\mathcal{A}$, we obtain that $\mathcal{A}$ is UGAS. Since the HDS is well-posed, by \cite[Lemma 7.20]{Goebel:12}  the UGAS property is robust and the result of item (b) holds.

To show inequality \eqref{sub_optimal_measure} in item (c), note that since the Lyapunov function \eqref{Lyapunov_function1} does not increase during flows and jumps, we have that $V(z(t+\ell,j))\leq V(z(t,j))$, for all $\ell\in[0,\Delta T]$, where $\Delta T=T_{\max}-T_{\min}$, and $(t,j)\in\text{dom}(z)$. Let $(t_d,j_d)\in\text{dom}(z)$ be the hybrid times such that $z(t_d,j_d)\in D$. Then, by construction $z(t_d+\Delta T,j_d+1)\in D$, and the Lyapunov function satisfies $V(z(t_d+\Delta T,j_d+1))\leq V(z(t_d,j_d+1))$, that is
\begin{align}
& V(t_d+\Delta T, j_d+1)\leq \frac{|x_2^{+}-x^*|^2}{2}+c(\tau^+)^2 [f(x_1^+)-f^*]\notag\\
&~=\frac{|x_1(t_d,j_d+1)-x^*|^2}{2}+cT_{\min}^2  [f(x_1(t_d,j_d+1))-f^*]\label{inequality_rate1}
\end{align}
where in the last equality we used the jump rule \eqref{jumps_periodic2}. By the definition of $(t_d,j_d)$, $\Delta T$, and the jump set $D$, we have that $\tau^2(t_d+\Delta T,j_d+1)=T_{\max}^2$. Using the definition of $V$, $j'_d:=j_d+1$, $\tilde{f}(t_d+\Delta T,j'_d):=f(x_1(t_d+\Delta T,j_d'))-f^*$, and the strong convexity of $f(\cdot)$, we get from \eqref{inequality_rate1} 
\begin{align}
\tilde{f}(t_d+\Delta T,j'_d)&\leq \frac{1}{T_{\max}^{2}}\Bigg[\frac{1}{c\mu}+T_{\min}^2\Bigg]\tilde{f}(x_1(t_d,j'_d),\label{decrease_cost_function}
\end{align}
where the last inequality follows by the strong convexity of $f(\cdot)$. Equation \eqref{decrease_cost_function} says that by the end of each period of flow after a jump, the cost function decreases by a constant factor of $k_0:=\frac{1}{T_{\max}^2}\left[\frac{1}{c\mu}+T_{\min}^2\right]$ which satisfies $0<k_0<1$ due to condition \eqref{inequality_period2}. Since $z(t_d,j_d)\in D$ implies $t_d=(j+1)\Delta T$ and $j_d=j$, for all $j\in\mathbb{Z}_{\geq0}$, using \eqref{decrease_cost_function} for all $(t_d,j_d)$ we have: 
\begin{equation*}
\tilde{f}(\Delta T+j\Delta T,j)\leq  k_0 \tilde{f}(j\Delta T,j),
\end{equation*}
 and since $\tilde{f}(j\Delta T,j)=\tilde{f}(j\Delta T,j-1)$, we get
\begin{align}\label{exponential_bound}
\tilde{f}(\Delta T+j\Delta T,j)&\leq  k^{j+1}_0 \tilde{f}(0,0).
\end{align}
For each fixed $j\geq1$, we have that during flows $\tilde{f}(t,j+1)\leq k_1\tilde{f}(t_{0,j+1},j+1)$, where $t_{0,j+1}$ is the smallest $t$ in the time domain such that $(t,j+1)\in\text{dom(z)}$, and where $k_1=\frac{1}{T_{\min}^{2}}\left[\frac{1}{c\mu}+T_{\min}^2\right]$. Using \eqref{exponential_bound} we get
\begin{equation}\label{error_boundthm2a}
\tilde{f}(t,j)\leq \tilde{f}(\Delta T+(j-1)\Delta T,j-1)\leq k_1k^{j}_0 \tilde{f}(0,0).
\end{equation}
Since the hybrid time domain of the system is periodic, for each $(t,j)\in\text{dom}(z)$ we have 
\begin{equation}\label{error_boundthm2b}
j\geq \frac{\max\{t+j-\Delta T,0\}}{\Delta T+1}=:\tilde{\alpha}(t+j).
\end{equation}
Let $\tilde{k}_0:=1-k_0$. Using \eqref{error_boundthm2a} and the fact that $f\in\mathcal{F}_{\mu,L}$, we get
%
%
%
%
\begin{align*}
\tilde{f}(x_1(t,j))&\leq k_1(1-\tilde{k}_0)^{j} \left[\tilde{f}(x_1(0,0))\right]\notag\\
&\leq 0.5 k_1L\exp\left(-\tilde{k}_0\tilde{\alpha}(t+j) \right) |\tilde{x}_1(0,0)|^2,
\end{align*}
where the last inequality follows by the Lipschitz gradient condition and where $\tilde{x}_1(0,0)=x_1(0,0)-x^*$. This inequality establishes item (c) with constants $k_a=0.5 k_1L$ and $k_b=\tilde{k}_0$. \null \hfill \null $\blacksquare$
%
%
%
%
\subsection{Proof of Lemma \ref{lemma_optimal_switching}}

Taking the derivative of $k^{\frac{t}{\Delta T}}_1$ with respect to $\Delta T$, and equating to zero, we obtain that $\Delta T^*=e~\sqrt{\frac{1}{c\mu}+T^2_{\min}}$, which satisfies \eqref{stronger_condition}.  Substituting $\Delta T^*$ in $k_1$, using the fact that $k_0\leq k_1$, and using $k_1$ instead of $k_0$ in the bound \eqref{error_boundthm2a}, we get $\tilde{f}(t,j)\leq k^{j}_1 \tilde{f}(0,0)$. For $\varepsilon>0$ we have that $k^j_1 \tilde{f}(0,0)<\varepsilon$ whenever $j\geq 0.5 \log\left(\frac{\tilde{f}(0,0)}{\varepsilon}\right)$. Multiplying by $\Delta T^*$ at both sides we get the result.  \null \hfill \null $\blacksquare$

\end{document}